\documentclass[11pt,dutch]{article}
\usepackage{amsmath}
\usepackage{amsthm}
\usepackage{amsfonts}
\usepackage{amssymb}
\usepackage[all]{xy}
\theoremstyle{plain}
\newtheorem{theorem}{Theorem}[section]
\newtheorem{lemma}[theorem]{Lemma}

\newtheorem{proposition}[theorem]{Proposition}

\theoremstyle{definition}
\newtheorem{definition}[theorem]{Definition}

\theoremstyle{remark}

\newtheorem{counterexample}[theorem]{Counterexample}

\newcommand{\monr}{\ar@{{^ (}->}[r]}
\newcommand{\monl}{\ar@{{^ (}->}[l]}
\newcommand{\monu}{\ar@{{^ (}->}[u]}
\newcommand{\mond}{\ar@{{^ (}->}[d]}
\newcommand{\monrr}{\ar@{{^ (}->}[rr]}
\newcommand{\monll}{\ar@{{^ (}->}[ll]}
\newcommand{\monuu}{\ar@{{^ (}->}[uu]}
\newcommand{\mondd}{\ar@{{^ (}->}[dd]}
\newcommand{\mondr}{\ar@{{^ (}->}[dr]}
\newcommand{\monur}{\ar@{{^ (}->}[ur]}
\renewcommand{\to}{\longrightarrow}

\newcommand{\del}{\partial}

\newcommand{\Bc}{\ensuremath{\mathcal{B}}}
\newcommand{\Ac}{\ensuremath{\mathcal{A}}}

\newcommand{\Pc}{\ensuremath{\mathcal{P}}}
\newcommand{\Ab}{\mathsf{Ab}}

\newcommand{\N}{\ensuremath{\mathbb{N}}}

\newcommand{\Gp}{\ensuremath{\mathsf{Gp}}}
\newcommand{\Rng}{\ensuremath{\mathsf{Rng}}}

\newcommand{\PXM}{\ensuremath{\mathsf{PXMod}}}
\newcommand{\XM}{\ensuremath{\mathsf{XMod}}}
\newcommand{\PX}{\ensuremath{\mathsf{PXM}}}
\newcommand{\X}{\ensuremath{\mathsf{XM}}}

\newcommand{\RG}{\ensuremath{\mathsf{RG}}}
\newcommand{\Gpd}{\ensuremath{\mathsf{Gpd}}}

\newcommand{\m}{\textit{\textbf{m}}}
\newcommand{\n}{\textit{\textbf{n}}}

\newcommand{\p}{\textit{\textbf{p}}}
\newcommand{\x}{\textit{\textbf{x}}}
\newcommand{\y}{\textit{\textbf{y}}}
\newcommand{\tv}{\textit{\textbf{t}}}
\newcommand{\rv}{\textit{\textbf{r}}}

\newcommand{\iv}{\textit{\textbf{1}}}
\newcommand{\av}{\textit{\textbf{a}}}
\newcommand{\bv}{\textit{\textbf{b}}}

\newcommand{\mv}{\textit{\textbf{m}}}
\newcommand{\nv}{\textit{\textbf{n}}}

\newcommand{\pv}{\textit{\textbf{p}}}
\newcommand{\xv}{\textit{\textbf{x}}}
\newcommand{\yv}{\textit{\textbf{y}}}
\newcommand{\zv}{\textit{\textbf{z}}}

\newbox\pullbackbox
\setbox\pullbackbox=\hbox{\xy \POS(4,0)\ar@{-} (0,0) \ar@{-} (4,4)
\endxy}

\newbox\pushoutbox
\setbox\pushoutbox=\hbox{\xy \POS(0,4)\ar@{-} (0,0) \ar@{-} (4,4)
\endxy}

\begin{document}

\CompileMatrices

\title{Relative Commutator Theory \\
in Varieties of $\Omega$-groups}
\author{T. Everaert\footnote{The author's research is financed by a Ph.D. grant of the institute of Promotion of Innovation through Science and Technology in Flanders.}}

\maketitle

\setcounter{section}{-1}
\begin{abstract}
\noindent We introduce a new notion of commutator which depends on a choice of subvariety in any variety of $\Omega$-groups. We prove that this notion encompasses Higgins's commutator, Fr\"ohlich's central extensions and the Peiffer commutator of precrossed modules.

\end{abstract}

\section{Introduction}

A  \emph{variety of $\Omega$-groups} \cite{Higgins} is a variety which has amongst its operations and identities those of the variety of groups but has no more than one constant. Examples are: the varieties of groups, (non unital) rings, commutative algebras, crossed modules and precrossed modules. In any variety of $\Omega$-groups a notion of commutator exists, introduced by Higgins \cite{Higgins}, which has as particular cases the ordinary commutators of groups and rings, amongst others.

In any variety $\Ac$, an algebra $A$  is called an \emph{abelian} algebra if it can be endowed with an internal group structure.  The subvariety of all abelian algebras in $\Ac$ will be denoted by $\Ac_{\Ab}$. If $\Ac$ is a variety of $\Omega$-groups then, for any $A\in\Ac$, $A$ is abelian if and only if the group operation on $A$ defines a homomorphism $A\times A\to A$. Higgins's commutator characterizes $\Ac_{\Ab}$ in the following way: for any $A\in\Ac$, the commutator $[A,A]$ is $\{1\}$ (the terminal algebra) if and only if $A\in\Ac_{\Ab}$.

An \emph{ideal} $N$ of an $\Omega$-group $A$ is a subalgebra of $A$ which is the kernel of some homomorphism $A\to B$. For such an ideal $N$, the quotient set $A/N$ admits a canonical $\Omega$-group structure. The corresponding inclusion $N\to A$ is called an \emph{extension}. In \cite{Froehlich} Fr\"ohlich introduced the following notion of central extension relative to a choice of subvariety $\Bc\leq\Ac$ in any variety of $\Omega$-groups $\Ac$: an extension of ${\Omega}$-groups $N\to A$ is a \emph{$\Bc$-central extension} if there is in $A$ only one element of the form $v(\nv\av)v(\av)^{-1}$, with $v\in W$, $\nv\in N$ and $\av\in A$, namely, the unit $1$. The ordinary central extensions of groups and of rings are particular cases of this notion. Furthermore, it is such that for any $\Omega$-group $A\in\Ac$, the identity $A\to A$ is a ($\Bc$-) central extension if and only if $A\in\Bc$. Furthermore, if $\Bc=\Ac_{\Ab}$ then an extension $N\to A$ is an $\Ac_{\Ab}$-central extension if and only if $[N,A]=\{1\}$.

In this paper we will define a new notion of commutator in any variety of $\Omega$-groups $\Ac$ relative to a choice of subvariety $\Bc\leq\Ac$ in such a way that it will characterize $\Bc$ (resp. the $\Bc$-central extensions) in the same way as Higgins's commutator characterizes $\Ac_{\Ab}$ (resp. the $\Ac_{\Ab}$-central extensions). More precisely, for any $A\in \Ac$, $A$ is in $\Bc$ if and only if the ($\Bc$-) commutator $[A,A]_{\Bc}$ is $\{1\}$; and, an extension $N\to A$ will be a $\Bc$-central extension if and only if $[N,A]_{\Bc}=\{1\}$. Higgins's commutator will be the particular case where $\Bc=\Ac_{\Ab}$.   

We will now give an example of a well known commutator in a variety of $\Omega$-groups which is \emph{not} a particular case of Higgins's commutator but will turn out to be a particular case of ours.

Recall that a \emph{precrossed} \emph{module} $(C,G,\del)$ is a group homomorphism $\del:C\to G$ equipped with a (left) group action of $G$ on $C$, such that
\[
\del(^gc) = g\del(c)g^{-1}
\]
for all $g\in G$ and $c\in C$. A morphism of precrossed modules $f:(C,G,\del)\to (D,H,\epsilon)$ is a pair of group homomorphisms $f_1:C\to D$ and $f_0:G\to H$ which preserve the action and are such that $\epsilon \circ f_1 =f_0\circ\del$. We write $\PXM$ for the category of precrossed modules and $\XM$ for the category of \emph{crossed} \emph{modules}, where this latter is the full subcategory of $\PXM$ whose objects $(C,G,\del)$ satisfy the condition
\[
^{\del(c)}c'=cc'c^{-1}
\]
for all $c, c'\in C$. This identity is often called the \emph{Peiffer identity}.

Also recall that a precrossed \emph{submodule} of a precrossed module $(C,G,\del)$ is a precrossed module  $(K,S,\kappa)$ such that $K$ and $S$ are, respectively, subgroups of $C$ and $G$, and such that  the action of $S$ on $K$ is a restriction of the action of $G$ on $C$ and $\kappa$ a restriction of $\del$ (in this case, we will erroneously write $\del$ instead of $\kappa$). $(K,S,\del)$ is a \emph{normal} precrossed submodule of $(C,G,\del)$ if, furthermore, $K$ and $S$ are, respectively,  normal subgroups of $C$ and $G$, and, for all $c\in C, g\in G, k\in K, s\in S$, one has $^{g}k\in K$ and $^{s}cc^{-1}\in K$. This is exactly the case when $(K,S,\del)$ is the kernel of some morphism $(C,G,\del)\to (D,H,\epsilon)$.

Let $(C,G,\del)$ be a precrossed module. The \emph{Peiffer commutator} of two normal precrossed submodules $(K,S,\del)$ and $(L,T,\del)$ of $(C,G,\del)$ is the normal subgroup of the group $K\vee L$, generated by the \emph{Peiffer elements} $\langle k,l\rangle =klk^{-1}(^{\del k}l)^{-1}$ and $\langle l,k\rangle =lkl^{-1}(^{\del l}k)^{-1}$, with $k\in K$ and $l\in L$. We will denote it by $\langle (K,S,\del),(L,T,\del)\rangle$.\footnote{Note that the Peiffer commutator is often defined as an ordinary subgroup of $C$. However, defining it as a \emph{normal} subgroup (of $K\vee L$) gives the commutator better properties and simplifies the comparison with other notions. Defining it as a normal subgroup of $K\vee L$ and not a normal sugroup of $C$ is natural, since defining it otherwise would yield the following ambiguous situation: suppose $(A',G',\del)$ is a precrossed submodule of $(A,G,\del)$, and  $(K,S,\del)$ and $(L,T,\del)$ normal precrossed submodules of both $(A,G,\del)$ and $(A',G',\del)$, then calculating the Peiffer commutator $\langle (K,S,\del),(L,T,\del)\rangle$ as a normal subgroup of $A$ or as a normal subgroup of $A'$ could give different results.}

Note that, although $\langle (K,S,\del),(L,T,\del)\rangle$ is usually considered a subgroup of $K\vee L$, it could as well be considered a precrossed submodule of $(K\vee L,S\vee T,\del)=(K,S,\del)\vee (L,T,\del)$; namely, the normal precrossed submodule $(\langle (K,S,\del),(L,T,\del)\rangle,\{1\},\del)$. We will adopt this point of view in Proposition \ref{Peiffer}.

It is well known that the category of precrossed modules is equivalent to a variety of $\Omega$-groups (see, e.g., \cite{LR} or \cite{Lo}). This allows us to compare the Peiffer commutator with Higgins's commutator in this particular variety. Since for any precrossed module $(C,G,\del)$, $\langle (C,G,\del),(C,G,\del)\rangle=\{1\}$ if and only if $(C,G,\del)$ is a crossed module and because $\PXM_{\Ab}\neq\XM$, the two commutators cannot be the same.  An explicit description of the abelian precrossed modules can be found, for example, in \cite{AL}. However, in Proposition \ref{Peiffer} we will prove that the Peiffer commutator is a particular case of the commutator introduced in this paper.

It is worth mentioning that the commutator of universal algebra introduced by Smith \cite{S} does \emph{not} coincide with Higgins's commutator and, consequently, is \emph{not} a particular case of the commutator introduced in this paper: although in any variety of $\Omega$-groups $\Ac$, for any $A\in\Ac$ and any ideal $N$ of $A$,
\[
[A,A] = \{1\} \Leftrightarrow [R_N,R_N] = \Delta_A \Leftrightarrow A\in\Ac_{\Ab},
\]
where $[A,A]$ denotes Higgins's commutator, $[R_N,R_N]$ Smith's commutator of the kernel equivalence relation of the quotient $A\to A/N$, and $\Delta_A$ the smallest equivalence relation on $A$, the implication
\[
[M,N] = \{1\} \Rightarrow [R_M,R_N] = \Delta_A
\]
does \emph{not} necessarily hold in $\Ac$ for all $A\in\Ac$ and all ideals $M$ and $N$ of $A$. A counterexample is given in \cite{Bourn-Commutator-Huq-Smith} in the variety of digroups. However, the induced notions of central extension do coincide \cite{MT}.

We would like to refer the reader interested in Fr\"ohlich's theory of central extensions and the related subject of \emph{Baer invariants} to the work of Lue \cite{Lue} and Furtado-Coelho \cite{Furtado-Coelho}, who further developed this theory. A categorical version of the notion of central extension was introduced by Janelidze and Kelly \cite{JK}, as an application of the categorical Galois theory developed by Janelidze \cite{J}. For a generalization of the theory of Baer invariants to the context of semi-abelian categories \cite{JMT}, we refer to Everaert and Van der Linden's papers \cite{EverVdL1} and \cite{EverVdL2}. Closely related to the present paper, in \cite{EG}, Everaert and Gran characterize the central extensions of internal precrossed modules over a fixed algebra $B$ in a semi-abelian variety $\Ac$, with respect to the subcategory of internal crossed modules in $\Ac$ over $B$. In the particular case where $\Ac$ is the variety of groups, this characterization is in terms of the Peiffer commutator. Already in \cite{Huq}, Huq introduced a categorical notion of commutator which, more recently, was generalized to regular Mal'tsev categories \cite{CLP} by Bourn in \cite{BCTRM}. In varieties of $\Omega$-groups, Higgins's commutator is easily seen to be equivalent to this commutator.

\section{Definition and Basic Properties}
For brevity, we will denote finite ordered sets $(x_1, x_2,\dots ,x_r)$, $(a_1, a_2,\dots ,a_s)$, $\dots$ by symbols $\xv, \av, \dots $ and write $f(\xv)$ for $f(x_1,x_2,\dots,x_r)$,  $f(\xv,\yv)$ for $f(x_1,\dots,x_r,y_1,\dots,y_s)$, etc. Instead of $(x_1y_1,x_2y_2,\dots,x_ry_r)$ we shall write $\xv\yv$.  Furthermore, $a_1,a_2,\dots,a_s\in A$ will be abbreviated to $\av\in A$. Also, if $1$ denotes the unit of a group operation, we shall write $\iv$ instead of $(1, 1, \dots, 1)$.

Throughout this paper, $\Ac$ will be a variety of ${\Omega}$-groups.

Whenever we use the word \emph{group} (resp. \emph{subgroup} or \emph{homomorphism}), it is understood that we mean \emph{$\Omega$-group} (resp. \emph{$\Omega$-subgroup} or \emph{homomorphism of $\Omega$-groups}), unless it is stated otherwise.

Now, let $\Bc$ be a subvariety of $\Ac$; then $\Bc$ is completely determined by a set of identities of terms. Since $\Ac$ is a variety of $\Omega$-groups, all these identities are of the form $v(\xv)=1$ and the corresponding terms $v(\xv)$ constitute a group \[
W=W_{\Bc}=\{ v(\xv) \ | \ v(\bv)=1,  \forall B\in \Bc, \forall \bv\in B\}
\]
(a subgroup of the group of terms $F_{\Ac}(\N)$, the free $\Omega$-group of $\Ac$ on a countable set). Of course, $B\in\Bc$ if and only if $v(\bv)=1$ for all $v\in W$ and $\bv\in B$.

Let us now introduce a new notion of commutator, which will be the object of study in this paper.

\begin{definition}\label{definitie}
For any $\Omega$-group $A\in\Ac$ and ideals $M, N$ of $A$, the commutator $[M,N]_{\Bc}$ is the ideal of $H=M\vee N$ generated by the set 
\[
\{ v(\mv \nv)v(\nv)^{-1}v(\mv)^{-1},  v(\p) \ | \ v\in W, \mv\in M, \nv\in N, \p\in M\wedge N\}.
\]
\end{definition}

We will sometimes abbreviate the term $v(\x \y)v(\y)^{-1}v(\x)^{-1}$ to $c_v(\x,\y)$ and write $C_{\Bc}(M,N)$ for the ideal of $M\vee N$ generated by the set 
\[
\{c_v(\m,\n) \ | \ v\in W, \m \in M, \n \in N\}.
\]

\begin{proposition}\label{eigenschappen}
For any $\Omega$-group $A\in\Ac$ and ideals $M, N, N'$ of $A$, one has:  
\begin{enumerate}
\item\label{reflectie}
$[A,A]_{\Bc}=\{1\} \Leftrightarrow A\in\Bc$;
\item\label{symmetrie}
$[M,N]_{\Bc} = [N,M]_{\Bc}$;
\item\label{doorsnede}
$[M,N]_{\Bc}\leq M\wedge N$;
\item\label{monotoon}
$N\leq N' \Rightarrow [M,N]_{\Bc}\leq [M,N']_{\Bc}$;
\item \label{universeel}universal property: 
if $q=q_{[M,N]_{\Bc}}: H= M\vee N\to H/[M,N]_{\Bc}$ is the canonical quotient, then $[q(M),q(N)]_{\Bc}=\{1\}$, where  $q(M)$ and $q(N)$ denote the direct images along $q$ of the subgroups $M$ and $N$, respectively; moreover,  $[M,N]_{\Bc}$ is the smallest ideal $I$ of $H$ such that $q_I$ has this property.  
\end{enumerate}
\end{proposition}
\begin{proof}
\ref{reflectie} follows readily from the definition.
In order to prove \ref{symmetrie}, we show that $v(\nv\mv)v(\mv)^{-1}v(\nv)^{-1}=c_v(\nv,\mv)\in [M,N]_{\Bc}$, for all $v\in W$, $\mv\in M$ and $\nv\in N$. We have:
\[
c_v(\nv,\mv)  =  c_v(\mv,\mv^{-1}\nv\mv) v(\mv)v(\mv^{-1}\nv\mv) v(\mv)^{-1}v(\nv)^{-1}
\]
\[
= c_v(\mv,\mv^{-1}\nv\mv) v(\mv) c_v(\mv^{-1}\nv\mv\nv^{-1},\nv) v(\mv^{-1}\nv\mv\nv^{-1})v(\nv)v(\mv)^{-1}v(\nv)^{-1}.
\]
Since $[M,N]_{\Bc}$ is an ideal of $H$, it suffices to observe that $c_v(\mv,\mv^{-1}\nv\mv)$, $c_v(\mv^{-1}\nv\mv\nv^{-1},\nv)$, $v(\mv^{-1}\nv\mv\nv^{-1}), v(\mv)v(\nv)v(\mv)^{-1}v(\nv)^{-1}\in [M,N]_{\Bc}$. In order to see that this last term is indeed in $[M,N]_{\Bc}$, let us abbreviate the term $v(\x)v(\y)v(\x)^{-1}v(\y)^{-1}$ to $w(\x,\y)$. Then, of course, $w\in W$, and, consequently,
\[
w(\mv,\nv) = w((\mv,\iv)(\iv,\nv)) = 
\]
\[
w((\mv,\iv)(\iv ,\nv)) w(\iv ,\nv)^{-1} w(\iv ,\mv)^{-1} =  c_w((\mv,\iv ),(\iv ,\nv)). 
\]
Hence, $w(\mv,\nv)$ is in $[M,N]_{\Bc}$.

Since $M$ and $N$ are ideals  of $H$, \ref{doorsnede} follows from the obvious identities 
\[
\frac{[M,N]_{\Bc}\vee M}{M}=\{1\}=\frac{[M,N]_{\Bc}\vee N}{N}. 
\]
\ref{monotoon} is an immediate consequence of the definition.

It follows from \ref{doorsnede} that $[q(M),q(N)]_{\Bc}=q[M,N]_{\Bc}=\{1\}$. Proving the second part of statement \ref{universeel} is straightforward.
\end{proof}

It is worth mentioning two negative results: the commutator defined above does \emph{not} preserve binary joins, i.e., in general,
\[
[M,N_1\vee N_2]_{\Bc}\neq [M,N_1]_{\Bc}\vee [M,N_2]_{\Bc};
\]
furthermore, the commutator is \emph{not} preserved by surjective images, i.e., in general,
\[
p([M,N]_{\Bc})\neq [p(M),p(N)]_{\Bc}
\]
(with $p:A\to B$ a surjective homomorphism). Let us give a counterexample to the first property:
\begin{counterexample}
Suppose $\Ac$ is the variety of commutative, not necessarily unital rings and $\Bc$ the subvariety determined by the identity $v_0(x)=x^3=0$; let $R$ be the free commutative ring generated by three elements $a_1$, $a_2$ and $b$, $R_1$ the ideal generated by $a_1$ and $b$, $R_2$ the ideal generated by $a_2$ and $b$, $S$ the ideal (of $R$) generated by $b$, and $I$ the ideal generated by  $a_1^2$, $a_2^2$ and $b^2$. Then, for every $r_1\in R_1$, $r_2\in R_2$ and $s\in S$, we have $v_0(s+r_1), v_0(r_1), v_0(r_2), v_0(s)\in I$ hence $v(s+r_1)-v(r_1)-v(s)\in I$ and $v(s+r_2)-v(r_2)-v(s)\in I$, for every $v\in W$. On the other hand, $v_0(b + a_1+a_2)-v_0(a_1+a_2)-v_0(b)\notin I$. Consequently 
\[
\Biggl[\frac{S\vee I}{I},\frac{R_1\vee I}{I}\Biggr]_{\Bc} = \{0\} = \Biggl[\frac{S\vee I}{I},\frac{R_2\vee I}{I}\Biggr]_{\Bc} 
\]
but
\[
\Biggl[ \frac{S\vee I}{I},\frac{R_1\vee R_2\vee I}{I}\Biggr]_{\Bc} \neq \{0\}.
\]
\end{counterexample}

Let us now take a closer look at the second (false) property. Suppose $\Bc$ is a subvariety of $\Ac$ such that the corresponding commutator \emph{is} preserved by surjective images. In this situation, it is readily seen that $\Bc\geq \Ac_{Ab}$. Indeed, suppose that the commutator is preserved by surjective images and $A\in\Ac_{Ab}$. In this case, the group operation on $A$ defines a (surjective) homomorphism $m:A\times A\to A$. Consequently, 
\[
[A,A]_{\Bc} = \Bigl[m \Bigl(A\times \{1\}\Bigr), m \Bigl(\{1\} \times A\Bigr)\Bigr]_{\Bc} = m \Bigl(\Bigl[A \times \{1\}, \{1\} \times A\Bigr]_{\Bc}\Bigr)
\]
\[
\leq m\Bigl(A \times \{1\} \wedge \{1\} \times A\Bigr) = \{1\}.
\]
It follows that $A\in \Bc$.

However, the converse is not true: the condition $\Bc\geq\Ac_{\Ab}$ does \emph{not} imply that the commutator corresponding to $\Bc$ is preserved by surjective images. In Proposition \ref{beeld} we will give a necessary and sufficient condition on the subvariety $\Bc$ for this property to hold. Before proving this proposition, we must recall some terminology from \cite{Higgins}. 
 
Suppose $t(\x,\y)$ is a term in two disjoint sets of indeterminates $\x$ and $\y$. It is called \emph{ideal term in $\x$ and $\y$},  if $t(\iv,\y)=1$.  For any subsets $M$ and $N$ of a group $A\in\Ac$, $M^N$ denotes the set 
\[
\{t(\m,\n) \ | \ t(\x,\y) \ \textrm{an ideal term in} \  \x  \  \textrm{and}  \ \y,  \ \m\in M, \n\in N\}.
\] 
$M^N$ is always an ideal of $M\vee N$ (where this latter denotes the subgroup of $A$ generated by $M\cup N$). In fact, it is the ideal of $M\vee N$ generated by $M$. Furthermore, if $p:A\to B$ is a surjective homomorphism, then $p(M^N)=p(M)^{p(N)}$. A term $t(\x,\y)$ in two disjoint sets of indeterminates $\x$ and $\y$ is called a \emph{commutator term in $\x$ and $\y$} if it is both an ideal word in $\x$ and $\y$, and in $\y$ and $\x$.

 \begin{lemma}\label{technisch}
 Let $F(X)$ and $F(Y)$ be free $\Omega$-groups in $\Ac$ on disjoint sets $X$ and $Y$. Then
 \[
 [ F(X)^{F(Y)}, F(Y)^{F(X)} ]_{\Bc} = C_{\Bc}(F(X)^{F(Y)}, F(Y)^{F(X)}).
 \]
 \end{lemma}
 \begin{proof}
 One inequality follows from the definition of the commutator. We prove the other one. 
 
Suppose $v\in W$ and $t_1(\x_1,\y_1), \dots, t_n(\x_n,\y_n) \in  F(X)^{F(Y)} \wedge F(Y)^{F(X)}=X^Y\wedge Y^X$ (with $\x_i \in X$ and $\y_i \in Y$ for all $i:1\dots n$); then each $t_i(\x_i,\y_i)$ is a commutator term in $\x_i$ and $\y_i$. Consequently, if we write $\tv(\x,\y)$ for $t_1(\x_1,\y_1), \dots, t_n(\x_n,\y_n)$, and $v_{\tv}(\x,\y)$ for $v(\tv(\x,\y))$, we have
 \[
 v(\tv(\x,\y)) = v\Bigl(\tv\Bigl( (\x,\iv)(\iv,\y) \Bigr)\Bigr) v\Bigl(\tv(\iv,\y)\Bigr)^{-1}  v\Bigl(\tv(\x,\iv)\Bigr)^{-1}
 \]
 \[
= c_{v_{\tv}}\Bigl( (\x,\iv),(\iv,\y) \Bigr)  \in  C_{\Bc}(F(X)^{F(Y)}, F(Y)^{F(X)}).\qedhere
 \]
 \end{proof}
 
 \begin{proposition}\label{beeld}
 For any subvariety $\Bc\leq \Ac$, the following conditions are equivalent:
 \begin{enumerate}
 \item
 for all surjective homomorphisms $p:A\to B$ in $\Ac$ and ideals $M$ and $N$ of $A$,
 \[
 p([M,N]_{\Bc}) = [p(M),p(N)]_{\Bc};
 \]
 \item
 for all $A\in \Ac$,  $C_{\Bc}(A,A)=  \{ v(\av) \ | \ v\in W, \av \in A\}$;
  \end{enumerate}
 \end{proposition}
\begin{proof}
Suppose 1 holds. In order to prove 2, it is sufficient to show that $[A,A]_{\Bc}= C_{\Bc}(A,A)$ for all $A\in\Ac$. Suppose $A\in\Ac$. Let us write $F(A)$ for the free group in $\Ac$ on the underlying set of $A$. Furthermore, let us denote by $F(A)\coprod F(A)$ the coproduct (the free product) of $F(A)$ with itself and by $\{1\} \coprod F(A)$ and $F(A) \coprod \{1\} $ the subgroups of $F(A)\coprod F(A)$ induced by the inclusion of $\{1\}$ into $F(A)$ and the identity on $F(A)$. From Lemma \ref{technisch} it follows that
\[
\Bigl[ \Bigl(F(A)\coprod \{1\} \Bigr)^{  ( \{1\} \coprod F(A) ) } \ , \ \Bigl(\{1\} \coprod F(A)\Bigr)^{  ( F(A) \coprod \{1\}  ) } \Bigr]_{\Bc}
\]
\[
= C_{\Bc}\Bigl( \Bigl(F(A)\coprod \{1\} \Bigr)^{  ( \{1\} \coprod F(A) ) } \ , \  \Bigl(\{1\} \coprod F(A)\Bigr)^{  ( F(A) \coprod \{1\} } \Bigr). 
\]
By assumption, the commutator is preserved by surjective images. Furthermore, it is readily seen that this is also the case for $C_{\Bc}(\cdot, \cdot)$.  Let us write $\epsilon_A$ for the unique homomorphism $F(A)\to A$ that sends each $a\in A$ to itself, and $(1_{F(A)},1_{F(A)})$ for the codiagonal $F(A)\coprod F(A)\to F(A)$. Then, applying $\epsilon_A\circ (1_{F(A)},1_{F(A)})$  to the identity above yields $[A,A]_{\Bc}= C_{\Bc}(A,A)$.

In order to see that 2 implies 1, it suffices to note that when 2 holds, the definition of the commutator simplifies to 
\[
[M,N]_{\Bc} = C_{\Bc}(M,N),
\]
for all $A\in \Ac$ and ideals $M, N$ of $A$. It has been remarked above that $C_{\Bc}(\cdot, \cdot)$ is preserved by surjective images. 
\end{proof}

We will now give an example of a subvariety $\Bc$ of a variety of $\Omega$-groups $\Ac$ that does not satisfy the conditions of Proposition \ref{beeld}, while $\Bc\geq \Ac_{Ab}$.

\begin{counterexample}
Suppose $\Ac$ is the variety of commutative, not necessarily unital rings and $\Bc$ the subvariety determined by the identity $v_0(x)=x^2=0$. It is well known that in this case $\Ac_{Ab}$ is determined by the identity $w_0(x,y)=xy=0$. Hence, one clearly has $\Bc\geq \Ac_{Ab}$. 

On the other hand, let $R$ be the free commutative ring generated by one element $a$. Obviously, $a^2\in [R,R]_{\Bc}$. We will show that it is not in $C_{\Bc}(R,R)$. Note that each term in $W$ is an additive sum of terms of the same form as $v_0(x)$ and terms of the same form as $v_t(x,\x)=x^2t(\x)$, for any term $t$. For any terms $t(\x)$ and $u(\y)$, let us write $t+u=(t+u)(\x,\y)$ for the term $t(\x)+u(\y)$. We have 
\[
c_{t+u}((\x,\y),(\x',\y')) = c_t(\x,\x') + c_u(\y,\y'),
\]
hence, $C_{\Bc}(R,R)$ is the ideal of $R$ generated by the sets 
\[
\{ c_{v_0}(r,r') \ | \ r, r'  \in R \}
\]
and 
\[
\{c_{v_t}( (r,\rv), (r',\rv') ) \ | \ t \ \textrm{a term,} \ r, r', \rv, \rv' \in R\}.
\] 
It is easily seen that the first of these sets is contained in the ideal of $R$ generated by $a^2+a^2$ and the second one in the ideal generated by $a^3$. Consequently, $a^2\notin C_{\Bc}(R,R)$. 

\end{counterexample}

\section{Particular Cases}

\subsection{Higgins's commutator}

In \cite{Higgins} Higgins defined a notion of commutator of $\Omega$-groups, which can be  characterized as follows (\cite{Higgins} Lemma 2.2(i) and Lemma 4.1): for any two ideals $M$ and $N$ of an ${\Omega}$-group $A$, the \emph{commutator} of $M$ and $N$, denoted by $[M,N]$, is the ideal of $H=M\vee N$ generated by  the set 
\[
\{ f(\mv\nv)f(\nv)^{-1}f(\mv)^{-1} \ | \ f \  \textrm{a term}, \mv\in M, \nv\in N\}.
\]
Since an $\Omega$-group $A$ is an abelian $\Omega$-group if and only if $[A,A]=\{1\}$, the subvariety $\Ac_{\Ab}$ is determined by the identities $f(\xv\yv)=f(\xv)f(\yv)$, for every term $f$. It is then readily seen that $[A,A]_{\Ac_{\Ab}}=[A,A]$.
 
We will now prove that this identity holds in general, for any two ideals $M$ and $N$ of an $\Omega$-group $A\in\Ac$.

\begin{proposition}\label{Higgins}
If $\Bc$ is the subvariety of all abelian groups, then the commutator defined in Definition \ref{definitie} is just Higgins's commutator of $\Omega$-groups. 
\end{proposition}
\begin{proof}
In order to prove $[M,N]_{\Ac_{\Ab}}\leq [M,N]$, we must prove that $v(\pv)\in [M,N]$, for all $v\in W$ and $\pv\in M\wedge N$. This follows from the monotony of Higgins's commutator. Indeed,  
\[
v(\pv)\in [M\wedge N,M\wedge N]_{\Ac_{\Ab}} = [M\wedge N,M\wedge N] \leq [M,N].
\]

In order to prove the other inclusion, we associate with every term $f(\xv)$ a term $c_f(\yv,\zv)=f(\yv\zv)f(\zv)^{-1}f(\yv)^{-1}$. Then $c_f\in W$ and, for all $\mv\in M$ and $\nv\in N$, we have: 
\[
f(\mv\nv)f(\nv)^{-1}f(\mv)^{-1}  = c_f(\mv,\nv) =c_f((\mv,\iv)(\iv,\nv))
\]
\[
=c_f((\mv,\iv)(\iv,\nv))c_f(\iv,\nv)^{-1}c_f(\mv,\iv)^{-1} \in [M,N]_{\Ac_{\Ab}}.\qedhere
\]
\end{proof}

As a corollary of the previous proposition, we get that the classical commutators of groups and rings, amongst others, are particular cases of the commutator introduced in this paper. 

Note that (the second part of) the proof of Proposition \ref{Higgins} implies also that $C_{\Bc}(A,A)=\{v(\av) | v \in W, \av \in A\}$, for every $A\in \Ac$. This is condition 2 of Proposition \ref{beeld}. Thus, as a corollary, we get that Higgins's commutator is preserved by surjective images. This, of course, comes as no surprise.

\subsection{Fr\"ohlich's central extensions}

As recalled in the introduction, Fr\"ohlich defined in \cite{Froehlich} a notion of central extension in any category of ${\Omega}$-groups $\Ac$, depending on a chosen subvariety of $\Ac$. In this subsection, we will show that this notion can be expressed in terms of our commutator, similar, e.g., to the situation in groups, where the central extensions can be characterized in terms of the ordinary group commutator.

 \begin{proposition}\label{Frohlich}
 An extension $N\to A$ in $\Ac$ is a $\Bc$-central extension in the sense of Fr\"ohlich if and only if $[N,A]_{\Bc}=\{1\}$. 
  \end{proposition}
 \begin{proof}
 Suppose $N\to A$ is a $\Bc$-central extension and suppose $v\in W$. We must prove that $v(\nv)=1=v(\nv\av)v(\av)^{-1}v(\nv)^{-1}$ for all $\nv\in N$ and $\av\in A$. Using the centrality twice, we first get $v(\nv)=v(\nv\iv)v(\iv)^{-1}=1$, and then $v(\nv\av)v(\av)^{-1}v(\nv)^{-1}=11=1$.
 
If, conversely, $[N,A]_{\Bc}=1$, then $v(\nv\av)v(\av)^{-1}=v(\nv\av)v(\av)^{-1}v(\nv)^{-1}=1$, for all $v\in W$, $\nv\in N$ and $\av\in A$, i.e. $N\to A$ is a $\Bc$-central extension.
  \end{proof}
 
 It is also readily seen that $[N,A]_{\Bc}$ is just Fr\"ohlich's \cite{Froehlich} $V_1$ of the extension $N\to A$ (or, using Furtado-Coelho's \cite{Furtado-Coelho} notation, $V_1(N | A)$).

\subsection{The Peiffer commutator}\label{PC}

As recalled in the introduction, the category of precrossed modules is equivalent to a variety of $\Omega$-groups (see, e.g., \cite{LR} or \cite{Lo}). In fact, it is equivalent to the variety whose theory consists, in addition to the group operation, the unit $1$ and the inversion $( \cdot )^{-1}$, of two $1$-ary operations $d$ and $c$, which satisfy the identities $d(xy)=d(x)d(y)$, $d(1)=1$,  $c(xy)=c(x)c(y)$, $c(1)=1$, $d\circ d=c\circ d=d$ and $d\circ c=c\circ c=c$. We will denote this variety by $\PX$. Remark that the identities imply, in particular, that the operations $d$ and $c$ are group homomorphisms $A\to A$, for any algebra $(A,d,c)\in\PX$.

With any $\Omega$-group $(A,d,c)\in\PX$ is associated the precrossed module $\Pc(A,d,c)=(K[d], I[d], c)$, where the action is given by
\[
^gk = gkg^{-1},
\]
for $k\in K[d]$, $g\in I[d]$. Conversely, with a precrossed module $(C,G,\del)$ is associated the $\Omega$-group $(G\ltimes C,d,c)$, where $d(g,c)= (g,1)$ and $c(g,c)=(g\del(c),1)$, and the product in $G\ltimes C$ is given by $(g,c)(g',c')=(gg',c{^{g}c'})$.

Via this equivalence $\Pc: \PX\to\PXM$, the category $\XM$ corresponds to the subvariety of $\PX$  of all $\Omega$-groups $(A,d,c)$ which satisfy $[K[d],K[c]]=1$, where this last commutator is the ordinary commutator of groups.

In fact, this equivalence is often presented as an equivalence between $\PXM$ and $\RG(\Gp)$, the category of reflexive graphs in $\Gp$ (see, e.g., \cite{LR} or \cite{Lo}). It is readily seen that $\PX$  is equivalent to $\RG(\Gp)$.

Similarly, the categories of precrossed rings and of crossed rings are equivalent to the category  $\RG(\Rng)$ of reflexive graphs in $\Rng$ and to the category $\Gpd(\Rng)$ of groupoids in $\Rng$, respectively \cite{LR}. Similar equivalences exist also in the case of commutative algebras, as discussed in \cite{Porter}.

We will now prove that the Peiffer commutator of precrossed modules is a particular case of the commutator defined in Definition \ref{definitie}; namely, the case where $\Ac$ is $\PX$ and $\Bc$ is $\X$. 

It is easily observed that, for any $(A,d,c)\in\PX$, $K[d] = \{ a^{-1}d(a)\ |\ a\in A\}$ and $K[c] = \{ a^{-1}c(a)\ |\ a\in A\}$. Consequently, in this case, the group $W$ is generated by all terms of the same form as
\[
v_0(x,y)=x^{-1}d(x)y^{-1}c(y)d(x)^{-1}xc(y)^{-1}y.
\]

It is also readily seen that, via the equivalence $\Pc$, precrossed submodules correspond to subgroups and normal precrossed modules to ideals.

\begin{proposition}\label{Peiffer}
If $\Ac$ is $\PX$ and $\Bc$ is $\X$ then the commutator defined in Definition \ref{definitie} corresponds to the Peiffer commutator via the equivalence $\Pc:\PX\simeq \PXM$. 
\end{proposition}
\begin{proof}
Suppose $(A,d,c)\in\PX$. We have to prove that, far any two ideals $(M,d,c)$ and $(N,d,c)$ of $(A,d,c)$, 
\[
\Pc([(M,d,c),(N,d,c)]_{\X})=\langle \Pc(M,d,c),\Pc(N,d,c)\rangle.
\]
By the universal properties of both commutators (see Proposition \ref{eigenschappen} (\ref{universeel})) it suffices to prove that 
\[
[(M,d,c),(N,d,c)]_{\X}=\{1\} \  \textrm{if and only if} \ \langle \Pc(M,d,c),\Pc(M,d,c)\rangle=\{1\}.
\]
Observe that for all $k\in K[d]\wedge M$ and $l\in K[d]\wedge N$,
\[
\langle k,l \rangle = klk^{-1}c(k)l^{-1}c(k)^{-1} = k[l,k^{-1}c(k)]k^{-1}
\]  
and, similarly, $\langle l,k \rangle = l[k,l^{-1}c(l)]l^{-1}$; here the square brackets denote ordinary group commutators. Since, moreover, 
\[
K[c]\wedge M = \{ k^{-1}c(k) \  | \ k\in K[d]\wedge M \} 
\]
and
\[
K[c]\wedge N = \{ l^{-1}c(l) \ | \ l\in K[d]\wedge N \},
\]
we have 
\[
\langle \Pc(M,d,c),\Pc(N,d,c)\rangle = \Bigl\langle (K[d]\wedge M,I[d]\wedge M,c), (K[d]\wedge N,I[d]\wedge N,c)\Bigr\rangle
\]
\[
 = \Bigl[K[d]\wedge N,K[c]\wedge M\Bigr] \vee \Bigl[ K[d]\wedge M,K[c]\wedge N\Bigr].
\]

Suppose that $[(M,d,c),(N,d,c)]_{\X}=1$. We prove that $[K[d]\wedge M,K[c]\wedge N]=1$. Observe that 
\[
K[d]\wedge M = \{m^{-1}d(m) \ | \ m\in M\}
\]
and
\[
K[c]\wedge N = \{n^{-1}c(n) \ |  \ n\in N\}.
\]
Consequently, if we prove for all $m\in M$ and $n\in N$ that
\[
v_0(m,n)=m^{-1}d(m)n^{-1}c(n)d(m)^{-1}mc(n)^{-1}n=1,
\]
the result follows. Since $v_0\in W$, this follows from the assumption.

Conversely, suppose $[K[d]\wedge N,K[c]\wedge M]\vee [K[d]\wedge M,K[c]\wedge N]=1$. Then, for every  $\mv\in M$ and $\nv\in N$,

\begin{eqnarray*}
v_0(\m\n) &  =& v_0(m_1n_1,m_2n_2)\\
&=&{n_1}^{-1}{m_1}^{-1}d({m_1})d({n_1}) {n_2}^{-1}{m_2}^{-1}c({m_2})c({n_2}) \\
&&d({n_1})^{-1}d({m_1})^{-1}{m_1}{n_1} c({n_2})^{-1}c({m_2})^{-1}{m_2}{n_2}\\
&=& {m_1}^{-1} d({m_1}) ( d({m_1})^{-1}{m_1}  {n_1}^{-1}{m_1}^{-1}d({m_1})d({n_1}) )  \\
&&( {n_2}^{-1}{m_2}^{-1}c({m_2}){n_2} )  {n_2}^{-1} c({n_2}) d({n_1})^{-1}d({m_1})^{-1}{m_1}\\
&& {n_1} c({n_2})^{-1}c({m_2})^{-1}{m_2}{n_2} \\
& =&  {m_1}^{-1} d({m_1}) ( {n_2}^{-1}{m_2}^{-1}c({m_2}) {n_2} )   \\
&&( d({m_1})^{-1}{m_1}  {n_1}^{-1}{m_1}^{-1}d({m_1})d({n_1})  ) {n_2}^{-1} c({n_2}) d({n_1})^{-1}\\
&&d({m_1})^{-1}{m_1}{n_1} c({n_2})^{-1}c({m_2})^{-1}{m_2}{n_2},
\end{eqnarray*}
since
\[
d({m_1})^{-1}{m_1}  {n_1}^{-1}{m_1}^{-1}d({m_1})d({n_1}) \in K[d]\wedge N 
\]
and
\[
 {n_2}^{-1}{m_2}^{-1}c({m_2}) {n_2} \in K[c] \wedge M.
\]
Then,
\begin{eqnarray*}
v_0(\m\n) & = &  {m_1}^{-1} d({m_1}) {m_2}^{-1} c({m_2}) ( c({m_2})^{-1} {m_2} {n_2}^{-1} {m_2}^{-1} c({m_2}) {n_2} )  \\
&&( d({m_1})^{-1}{m_1} )  {n_1}^{-1}{m_1}^{-1}d({m_1})d({n_1}) {n_2}^{-1} c({n_2}) d({n_1})^{-1}\\
&&d({m_1})^{-1}{m_1} {n_1} c({n_2})^{-1}c({m_2})^{-1}{m_2}{n_2}\\
&=&  {m_1}^{-1} d({m_1}) {m_2}^{-1} c({m_2}) ( d({m_1})^{-1}{m_1} ) \\
&&( c({m_2})^{-1} {m_2} {n_2}^{-1} {m_2}^{-1} c({m_2}) {n_2} ) {n_1}^{-1}{m_1}^{-1}d({m_1})d({n_1}) \\
&&{n_2}^{-1} c({n_2})  d({n_1})^{-1}d({m_1})^{-1}{m_1}{n_1} c({n_2})^{-1}c({m_2})^{-1}{m_2}{n_2},
\end{eqnarray*}
since
\[
c({m_2})^{-1} {m_2} {n_2}^{-1} {m_2}^{-1} c({m_2}) {n_2} \in K[c]\wedge N \quad {\rm and} \quad d({m_1})^{-1}{m_1} \in K[d]\wedge M.
\]
Finally, since
\[
{m_1}^{-1}d({m_1}) \in K[d]\wedge M \quad {\rm and} \quad d({n_1}){n_2}^{-1} c({n_2})d({n_1})^{-1} \in K[c]\wedge N,
\]
and both 
\[
{m_2}^{-1} c({m_2}) \in K[c]\wedge M 
\]
and
\[
{n_2}  {n_1}^{-1}  d({n_1}) {n_2}^{-1} c({n_2}) d({n_1})^{-1} {n_1} c({n_2})^{-1} \in K[d]\wedge N,
\]
we have
\begin{eqnarray*}
v_0(\m\n) & = &  {m_1}^{-1} d({m_1}) {m_2}^{-1} c({m_2})  d({m_1})^{-1}{m_1}   c({m_2})^{-1} {m_2} {n_2}^{-1} \\
&&{m_2}^{-1} c({m_2}) {n_2}  {n_1}^{-1}  ( {m_1}^{-1}d({m_1}) ) ( d({n_1}) {n_2}^{-1} c({n_2}) \\
&&d({n_1})^{-1} )  ( d({m_1})^{-1}{m_1} )  {n_1} c({n_2})^{-1}c({m_2})^{-1}{m_2}{n_2}\\
&=&  {m_1}^{-1} d({m_1}) {m_2}^{-1} c({m_2})  d({m_1})^{-1}{m_1}   c({m_2})^{-1} {m_2} {n_2}^{-1}\\
&& {m_2}^{-1} c({m_2}) {n_2}  {n_1}^{-1}   ( d({n_1}) {n_2}^{-1} c({n_2}) d({n_1})^{-1} ) \\
&& {n_1} c({n_2})^{-1}c({m_2})^{-1}{m_2}{n_2}\\
&=&  {m_1}^{-1} d({m_1}) {m_2}^{-1} c({m_2})  d({m_1})^{-1}{m_1}   c({m_2})^{-1} {m_2} {n_2}^{-1}\\
&& ( {m_2}^{-1} c({m_2}) )  ( {n_2}  {n_1}^{-1}  d({n_1}) {n_2}^{-1} c({n_2}) d({n_1})^{-1} {n_1} c({n_2})^{-1}  ) \\
&&  ( c({m_2})^{-1}{m_2} )  {n_2}\\
&=&  {m_1}^{-1} d({m_1}) {m_2}^{-1} c({m_2})  d({m_1})^{-1}{m_1}   c({m_2})^{-1} {m_2}  \\
&&{n_2}^{-1}( {n_2}  {n_1}^{-1}  d({n_1}) {n_2}^{-1} c({n_2}) d({n_1})^{-1} {n_1} c({n_2})^{-1}  )  {n_2}\\
&=& v_0(\m)v_0(\n)
\end{eqnarray*}
Similarly, one proves that $v_0(\m\n)=v_0(\n)v_0(\m)$. Consequently, we also have $v_0(\m)v_0(\n) =v_0(\n)v_0(\m)$.

Since $d(v_0(\xv))=c(v_0(\xv))=1$, $W$ is just the (ordinary) group generated by all terms of the same form as $v_0$. Hence, the identity $v(\mv\nv)=v(\mv)v(\nv)$ follows from the identities above, for \emph{all} $v\in W$ (and all $\mv\in M$ and $\nv\in N$). It remains to be shown that $v(\pv)=1$, for all $\pv\in M\wedge N$. Again, it suffices to prove this identity in the case that $v=v_0$,  but in this case it is clear. 
\end{proof}

Note that combining the previous proposition with Proposition \ref{Frohlich} yields the following result:  an extension $(K,S,\del)\to (A,G,\del)$ of precrossed modules is an $\XM$-central extension if and only if the Peiffer commutator $\langle (K,S,\del), (A,G,\del)\rangle$ is $\{1\}$. For precrossed modules over a fixed group $G$, a similar result was obtained in \cite{EG}, as a particular case of a more general theorem. 

Note also that 
\[
v_0(x,y)= v_0( (x,1)(1,y))v_0(1,y)^{-1}v_0(x,1)^{-1} = c_{v_0}((x,1),(1,y))
\]
implies that condition 2 of Proposition \ref{beeld} is satisfied in the case where $\Ac=\PX$ and $\Bc=\X$. As a corollary, we find that the Peiffer commutator is preserved by surjective images. This, of course, is nothing new.

\section*{Acknowledgements}
I would like to thank Tim Van der Linden for many helpful remarks and George Janelidze for pointing me to the fact that the commutator defined in this text is not preserved by surjective images.

\noindent Tomas Everaert \\
{ Vakgroep Wiskunde \\
 Faculteit Wetenschappen \\
Vrije Universiteit Brussel \\
Pleinlaan 2, 1050 Brussel \\
Belgium \\}
Email: \texttt{teveraer@vub.ac.be}

\end{document}